\documentstyle{amsppt}
\NoBlackBoxes
\leftheadtext{Pei-Kee LIN}
\rightheadtext{UNRESTRICTED PRODUCTS OF CONTRACTIONS} 
\topmatter 
\title UNRESTRICTED PRODUCTS OF CONTRACTIONS
IN BANACH SPACES  \endtitle 
\author 
Pei-Kee Lin \endauthor
\address Department of
Mathematics, Memphis State University, Memphis, TN  38152      
\endaddress \email LINPK\@hermes.msci.memst.edu \endemail

%  Math Subject Classifications
\subjclass  47A05, 47B05, 65J10 \endsubjclass\keywords Contraction,
weak convergence, random product \endkeywords
\abstract 
Let $X$ be a reflexive Banach space such that for any $x \ne 0$ the
set 
$$
\{x^* \in X^*: \text {$\|x^*\|=1$ and $x^*(x)=\|x\|$}\}
$$
is compact.
We prove that any unrestricted product of of a finite number of
$(W)$ contractions on $X$ converges weakly.
\endabstract     
\endtopmatter
     
\document

\def \AA {1}
\def \Ba {2}
\def \Bb {3}
\def \DKR {4}
\def \DRa {5}
\def \DRb {6}
\def \DRc {7}
\def \vN   {8}
  Recall a (bounded) linear operator on  a  Banach space
$X$ is said to be {\it contraction} if $\|Tx \| \leq \|x\|$ for
all $x \in X$.  $T$ is said to satisfy {\it condition} $(W$) if
$\{x_n\}$ is bounded and $\|x_n\|-\|T x_n\|$ converges
to 0 implies $x_n -T x_n$ converges to 0 weakly.  An algebraic
semigroup ${\bold S}$ generated by a (possibly infinite) set of
contractions is said to satisfy {\it condition} $(W)$ if for any
bounded sequence of vectors $\{x_n\}\subseteq X$ and a sequence of
words $\{W_n\}$ from ${\bold S}$ such that $\|x_n\|-\|W_n x_n\|$
converges to $0$, $x_n-W_n x_n $ converges to 0 weakly.
Let
$\{ T_1, T_2, \ldots , T_N\}$ be $N$
$(W)$-contractions on $X$,
and let
$r$ be a mapping from the set of natural numbers ${\Bbb N}$ onto
$\{1, 2, \ldots , N\}$, which assumes each value infinitely often.
An {\it unrestricted} (or {\it random}) {\it product} 
of these operators is the
sequence $\{ S_n:\enspace n = 1, 2, \ldots \}$ defined
by 
$$
S_n = T_{r(n)} T_{r(n-1)} \ldots T_{r(1)}.
$$ 
John von Neumann [\vN] proved that if $T_1$ and $T_2$ are
orthogonal projections on Hilbert space, and if $\{S_n\}$ is a
random product of $\{T_1, T_2\}$, then $S_n$
converges strongly.   Amemiya and Ando [\AA] extended this result
by showing if $\{S_n\}$ is a random product of finite
$(W)$-contractions
on Hilbert space, then  $S_n   $ converges weakly.   
On the other hand, Bruck [\Ba] showed that a random product of
infinite $(W)$-contractions on Hilbert space is not necessary to
be convergent weakly.        One may ask 
the following question.
\example{Question 1}  Does every random product of finite number of
$(W)$-contractions on a reflexive Banach space $X$ converge weakly?
\endexample

Recently, J. M. Dye, A. Khamsi, and S.
Reich [\DKR] showed the answer is positive if $X$ is smooth.
Indeed, they proved that if $X$ is a smooth reflexive Banach
space, and $\{T_1,T_2, \cdots,T_N\}$ are $N$ $(W)$-contractions
on $X$, then 
\roster
\item "(i)" there is a contraction projection $Q$ of $X$ onto
the common fixed point set such that $QT_j=T_j Q$ for $1 \leq j
\leq N$;
\item "(ii)" the semigroup ${\bold S}=S(T_1,T_2,\cdots,T_N)$
generated by $\{T_1,T_2,\cdots,T_N\}$
has property $(W)$;
\item "(iii)" for any random
product $S_n$ drawn from $\{T_1,T_2,\cdots,T_N\}$, 
$S_n x$ converges to $Qx$ weakly.
\endroster

But the following example shows (i) is not true in general.

\example{Example 1} Let $X=\ell_1^2$ and let $\{e_1,e_2\}$ be the
natural basis of $\ell_1^2$.  $T_1$ and $T_2$ are the
contractions defined by
$$
\align
T_1(a e_1 +b e_2)=&(a  + \frac b 2) e_1 \\
T_2(a e_1 +b e_2)=&(a  + \frac b 3) e_1.
\endalign
$$
It is easy to see that $T_1$ and $T_2$ are $(W)$-contractions.
But both $T_1$ and $T_2$ are contraction porjections onto
the common fixed point set.  Hence, there is no contraction projection 
$Q$ onto
$\{a e_1: a \in \Bbb R\}$ such that $Q$ commutes with $T_1$ and 
$T_2$. But any random product drawn from $T_1$ and $T_2$
converges weakly (strongly).
\endexample
In this article, we modify their ideas and we
prove that
if 
$X$ is a reflexive Banach space such that for any $x \ne 0$ the
set 
$$
\{x^* \in X^*: \text {$\|x^*\|=1$ and $x^*(x)=\|x\|$}\}
$$
is compact,
then any random product of finite number of $(W)$-contractions on
$X$ converges weakly.

The proof is base on the following four lemmas.

\proclaim{Lemma 1}  Let $\{T_1, T_2,\cdots,T_N\}$ be $N$
contractions and let $Y$ be the common fixed point set.  For
any bounded sequence $\{x_n:n \in \Bbb N\}$ in $X$ and
any sequence of words $\{W_n:n \in \Bbb
N\}$  from the semigroup 
${\bold S}=S(T_1,T_2,\cdots,T_N)$,
let 
$Z^*$ be the set consisting of all $x^*$ such that
\roster
\item "(iv)"  $x^*$ is supporting at some point $y \in Y$;
\item "(v)" the set $\{W_n^* x^*: n \in \Bbb N\}$ is relatively 
compact.
\endroster
If $Z^*$ separates the points of $Y$ and if both $\{x_n\}$ and $\{W_n
x_n\}$ converge weakly to some points in $Y$, say $u$ and $v$.
Then
$u=v$.
\endproclaim
\demo{Proof} Suppose that $u \ne v$.  Then there is $z^* \in Z^*$
such that $z^*(u) \ne z^*(v)$.  Since $\{W^*_n z^*: n \in \Bbb
N\}$ is compact, by passing to a subsequence if necessarily, we
may assume that $W^*_n z^*$ converges to $w^*$.  Since $Y$ is the
common fixed point set, 
$$W^*_n z^* |_Y=z^*|_Y, \quad \text {and} \quad
z^*|_Y=w^*|_Y.$$
Hence,
$$
z^*(v)=\lim_{n\to \infty} z^*(W_n x_n)=\lim_{n\to \infty}(W^*_n
z^*)(x_n)=\lim_{n \to \infty} w^*(x_n)=w^*(u)=z^*(u).
$$
This is a contradiction.  So $u=v$. \qed
\enddemo

\proclaim{Lemma 2} Let $T$ be a (W)-contraction on $X$ and let
$\{x_n\}$ be a  weakly convergent sequence.  If $\|x_n\|-\|T
x_n\|$ converges to 0, then the weak limit of $\{x_n\}$ is a
fixed point of $T$.
\endproclaim
\demo{Proof} Since $T$ is linear,
$$
\text {w-}\lim_{n \to \infty} T x_n=T(\text {w-}\lim_{n \to
\infty} x_n).
$$
On the other hand, $T$ is a $(W)$-contraction and $\|x_n\|-\|T
x_n\|$ converges to 0. So
$$
\text{w-}\lim_{n \to \infty} (x_n-T x_n)=0.
$$
This implies w-$\lim_{n \to \infty} x_n$ is a fixed point of $T$.
\qed \enddemo

\proclaim{Lemma 3} $($Proposition 6 $[\DKR])$  Let $\{T_1, T_2,
\cdots, T_N\}$ be $N$ $(W)$-contractions such that for any
proper subset $A$ of $\{1,2,\cdots,N\}$, the semigroup generated
by $\{T_j:j \in A\}$ satisfies property $(W)$.  Suppose that
$\{x_n\}$ is a bounded sequence in $X$ and $\{W_n\}$ is a
sequence of words from ${\bold S}=S(T_1,T_2,\cdots,T_N\}$ such
that
\roster
\item "(vi)" $\{x_n\}$ converges weakly to $u$, and $\{W_n x_n\}$
converges weakly to $v$;
\item "(vii)" $\lim_{n \to \infty} \|x_n\|-\|W_n x_n\|=0$.
\endroster
If $u \ne v$, then both $u$ and $v$ are common fixed points
of $\{T_1,T_2,\cdots, T_N\}$.
\endproclaim
\demo{Proof}  We claim that  $W_k$ is
complete (in $T_1,T_2,\cdots,T_N$) except finite number of $k$.

Assume it is not true.  By passing to a subsequence of $\{W_n\}$
and reindexing of the $T_j$'s,
we may assume that  $W_n \subseteq S(T_1,T_2,\cdots,T_{N-1})$. 
But $S(T_1,T_2,\cdots,T_{N-1})$ has property $(W)$.
(vii) implies
$$
\text {w-} \lim_{n \to \infty} (x_n-W_n x_n)=0.
$$
Hence, $u=v$.  This is a contradiction.

By passing to a subsequence of $\{W_n\}$, we may assume there
exist  words $T_j$, $T_{j'}$, $F_n  \in S(T_1,\cdots,T_{j-1},
T_{j+1}, \cdots,T_N)$, $F_n'  \in S(T_1,\cdots,T_{j'-1},
T_{j'+1}, \cdots,T_N)$, and $W'_n \in {\bold
S}=S(T_1,\cdots,T_N)$ such that
$W_n=(\cdots)T_j F_n =F'_n T_{j'} W'_n$.

Note:
$$
\align
\|x_n\|-\|W_n x_n\| \geq & \|x_n\|-\|T_j F_n x_n\|
 \geq \|x_n \|-
\|F_n x_n\|,\\
\|x_n\|-\|W_n x_n\| \geq & \|x_n\|-\|T_j F_n x_n\|
 \geq \|F_n x_n \|-
\|T_j F_n x_n\|,\\
\|x_n\|-\|W_n x_n\| \geq & \|W_n' x_n\|-\|W_n x_n\|
\geq \|T_{j'}
W_n' x_n\|-\|W_n x_n\|,\\
\|x_n\|-\|W_n x_n\| \geq & \|W_n' x_n\|-\|W_n x_n\|
\geq \|
W_n' x_n\|-\|T_{j'} W'_n x_n\|.
\endalign
$$
But the semigroups $S(T_1,\cdots,T_{j-1},
T_{j+1}, \cdots,T_N)$, and $ S(T_1,\cdots,T_{j'-1},
T_{j'+1}, \cdots,T_N)$ have property $(W)$. So
$$
\align
u=\text {w-}\lim_{n \to \infty} x_n=& \text {w-}\lim_{n \to \infty}
F_n x_n
= \text {w-}\lim_{n \to \infty}
T_j F_n x_n\\
\text {w-}\lim_{n \to \infty}W_n' x_n= &\text {w-}\lim_{n \to \infty}
T_{j'}W'_n x_n
= \text {w-}\lim_{n \to \infty}
F'_n T_{j'} W'_n x_n=v.
\endalign
$$
By Lemma 2,  $u$ (respectively, $v$) is a fixed point of $T_j$
(respectively, $T_{j'}$).  We claim that $u$ is a common fixed
point.

For each $k \ne j$, there exist words $G_n$ such that
$F_n=(\cdots) T_k G_n$. Since the semigroup $S(T_1,\cdots,T_{j-1},
T_{j+1}, \cdots,T_N)$ has property $(W)$ and
$\lim_{n \to \infty} \|G_n x_n\|-\|T_k G_n x_n\|=0$,
$u$ is a fixed point of $T_k$.   This implies $u$ is a common fixed point.
Similarly, $v$ is a common fixed point. \qed \enddemo

The following lemma was proved by J. Dye, M. A. Khamsi and S.
Reich.

\proclaim {Lemma 4} 
%$(1)$ $($Proposition 6 $[\DKR])$ 
%If $T$ is a $(W)$-contraction on a reflexive $X$,
%then the semigroup generated by $T$ satisfies property $(W)$.
%
$($Theorem 2 $[\DKR])$ Let ${\bold
S}=S(T_1,T_2,\cdots,T_N)$ be the semigroup generated by $N$
contractions on a reflexive Banach space $X$.
If ${\bold S} $ has property $(W)$, then every random product
$S_n$ drawn from $\{T_1,T_2,\cdots,T_N\}$ converges weakly.
\endproclaim
\demo{Proof} 
%Proof of (1). By the mean ergodic theorem, there is
%a contraction projection $Q$ from $X$ onto the fixed point set of
%$T$ such that $QT=Q=TQ$.  Suppose $\{x_k\}$ is a bounded sequence
%in $X$ and $\{n_k\}$ is a sequence of $\Bbb N$ such that
%$\|x_k\|-\|T^{n_k} x_k\|$ converges to 0.  We claim that $x_k-
%T^{n_k} x_k$ converges to 0 weakly.
%
%Suppose it is not true.  By passing to a subsequence if
%necessary, we may assume that $x_k$ converges to $u$ weakly,
%$T^{n_k} x_k$ converges to $v$ weakly, and $u \ne v$.  By
%Lemma 2, both $u$ and $v$ are fixed point of $T$.
%
%Subclaim.  If $\{y_k\}$ converges to a fixed point of $T$, then
%$\{(I-Q) y_k\}$ converges to 0 weakly.
%
%Assume that Subclaim were proved. Then
%$$
%u=\text {w-}\lim_{k \to \infty} Q x_k =\text {w-}\lim_{k \to
%\infty}T^{n_k}\circ  Q x_k =\text {w-}\lim_{k \to
%\infty} Q\circ T^{n_k} x_k =v.
%$$
%This is a contradiction.  So $u=v$.
%
%Suppose Subclaim is not true.  By passing to a subsequence of
%$\{y_k\}$ if necessary, we may assume that w-$\lim_{k \to \infty}
%(I-Q) y_k=w' \ne 0$.  Let $w=\text {w-}\lim_{k \to \infty} y_k
%\in Q(X)$.  Then
%$$
%w-w'=\text {w-}\lim_{k \to \infty} Q y_k \in Q(X).
%$$
%This is impossible.  And we proved Subclaim.
%
Suppose it is not true.  There exist a random
product $S_n$, $x \in X$ and two subsequences $\{n_k\}$, $\{m_k\}$ of $\Bbb
N$ such that $$\text {w-}\lim_{k \to \infty} S_{n_k}(x)=u\ne v =\text
{w-}\lim_{k \to \infty} S_{m_k}(x).$$
By passing to a subsequence of  $\{m_k\}$ if
necessary, we may assume that for every $k \in \Bbb N$, $n_k \leq
m_k$.  Let $W_k\in {\bold S}$ be the words such that $S_{m_k}=W_k
S_{n_k}$. Note: $\{\|S_n x\|\}$ is a decreasing sequence.  So 
$$ \lim_{k
\to \infty} \|S_{n_k} x\|-\|S_{m_k} x\|=0.
$$
But ${\bold S}$ has property $(W)$.
$$
u-v=\text {w-}\lim_{k \to \infty} S_{n_k} x-S_{m_k} x
=\text {w-}\lim_{k \to \infty} S_{n_k} x-W_kS_{n_k} x=0.
$$
This is a contradiction.  So $\{S_n x\}$ converges weakly.
\qed \enddemo

\proclaim {Theorem 5} Let $X$ be a reflexive Banach space such that for any $x \ne 0$ the
set 
$$
\{x^* \in X^*: \text {$\|x^*\|=1$ and $x^*(x)=\|x\|$}\}
$$
is compact.
If ${\bold S}$ be the semigroup generated by $N$ 
$(W)$-contractions $\{T_1,T_2,\cdots,T_N\}$ 
on $X$, then ${\bold S}$ has property $(W)$.
Hence, every random product drawn from $\{T_1,T_2,\cdots,T_N\}$
converges weakly.
\endproclaim
\demo{Proof} Suppose the Theorem is not true.  Let $N$ be the
smallest number such that there exists a semigroup $\bold S$
generated by $N$ $(W)$-contractions $\{T_1,T-2,\cdots,T_N\}$ 
without property $(W)$.  We
note that the semigroup generated by empty set is $\{I\}$ and it
has property $(W)$.

Let $\{x_n\}$ be a bouned sequence in $X$ and $W_n \in {\bold S}$
be a sequence of words such that
$$
\align
&\lim_{n \to \infty} \|x_n\|-\|W_n x_n\|=0\\
&\text {w-}\lim_{n\to \infty} x_n-W_n x_n \ne 0.
\endalign
$$
By passing to a subsequence if necessary, we may assume that
$$
\text {w-}\lim_{n \to \infty} x_n=u \ne v =\text {w-} \lim_{n \to
\infty} W_n x_n.
$$
By Lemma 3, both $u$ and $v$ are common fixed points.  Let $Y$ be
the common fixed point set.  Let  $Z^*$ be the set
$$
\align
\{x^*\in X^*:& |x^*(y)|=\|x^*\|\|y\| \text { for some $y \in Y
\setminus 0$ and}\\
&\text { $\{W^*_n x^*\}$ contains a convergent
subsequence} \}.
\endalign$$
For any non-zero common fixed point $y$, if $z^*$ is a 
a support functional of $y$, then $W^* z^*$ is a support
functional of $y$.  Assumption implies
$Z^*$ separates the points
of $Y$.  By Lemma 1, we have $u=v$. This is a contradiction.  Hence,
$\bold S$ has property $(W)$. \qed \enddemo

Recall a contraction $T$ is said to satisfy {\it condition}
$(W')$ if $\|x\|=\|T x\|$ implies $x=Tx$.  Suppose $T$ is a 
$(W)$-contraction such that $T^*$ satisfies property $(W')$.  If
$x^*$ is a supporting functional of some fixed point $x$, then
$T^* x^*$ is also a supporting functional of $x$.  But $T^*$ has
property $(W')$. This implies $T^*x^*=x^*$.  Hence, we have the following
Theorem.

\proclaim {Theorem 6} $($Remark of Theorem 1 $[\DKR])$
Let $\{T_1,T_2,\cdots,T_N\}$ be $N$ $(W)$-contractions on a
reflexive Banach space such that each $T_j^*$ has property
$(W')$.  Then the semigroup generated by $T_1,T_2,\cdots,T_N$ has
property $(W)$.
\endproclaim

\remark {Remark 1} Indeed, under the assumption of Theorem,
Theorem 1 and its Remark of [\DKR] show there is a contraction
projection $Q$ which commutes with $T_j$ for $1 \leq j \leq N$.
\endremark

One may ask when the random product converges strongly.  The
following lemma shows that it is enough to prove it contains a
convergent subsequence.

\proclaim{Lemma 7}  Let $\{T_1,T_2,\cdots,T_N\}$ be $N$ 
$(W')$-contractions, and let $S_n$ be a random product drawn from
them.  If there is a subsequence of $\{S_n x\}$ converges strongly,
then $\{S_n x\}$ converges strongly.
\endproclaim
\demo{Proof} Suppose that $\{S_{n_k} x\}$ converges to $z$.  We
claim that $z$ is a common fixed point.

Suppose the claim were proved.  Then $\|S_n x - z\|$ is a
decreasing sequence. So 
$$ \|S_n x - z\|= \|S_{n_k} x -z\|=0.
$$
This implies $S_n x$ converges to $z$.

Suppose $z$ is not a common fixed point.  By passing to a
subsequence of $\{n_k\}$ if necessary, we may assume that there
exist a sequence $\{m_k\}$, words $W_k$, $T_j$ from $\bold S$ such that
\roster
\item "(viii)"  $n_k \leq m_k \leq n_{k+1}$;
\item "(ix)" $S_{m_k+1} =T_j S_{m_k}=T_j W_k S_{n_k}$;
\item "(x)" $W_k z=z$ and $T_j z \ne z$.
\endroster
Since $W_k$'s are contractions, $\{S_{m_k} x\}$ converges to $z$.
If $z \ne T_j z$, then $\|T_j z\| <\|z\|$.  This implies
$$
\lim_{k \to \infty} \|S_{n_k} x\|=\|z\| > \|T_j z\|=\lim_{k \to
\infty}\| S_{m_k +1} x \|.
$$
This contradicts that $\{\|S_n x\|\}$ is a convergent sequence.
So $z$ must be a common fixed point. \qed
\enddemo

Hence, we have the following theorem.

\proclaim{Theorem 8}  Let $\{T_1,T_2,\cdots,T_N\}$ be $N$ 
$(W')$-contractions.  If $T_1$ is compact, 
then any random product $S_n$ drawn from
them
converges strongly.
\endproclaim

\remark{Remark 2} The above Theorem is still true if
$\{T_1,T_2,\cdots,T_N\}$ are $N$ $(W')$-nonexpansive mappings.
For more results of nonlinear nonexpansive mappings, see [\DRb],
[\DRc], and their reference.
\endremark


\Refs
    \ref \key \bf{\AA}  \by I. Amemiya and T. Ando \paper Convergence of 
random products of contractions in Hilbert space 
\jour Acta Sci. Math. (Szeged) \vol 26 \yr 1965 \pages239--244 \endref\ref \key \bf {\Ba} \by R. E. Bruck \paper Random products of
contractions in matric  and Banach spaces \jour J. Math. Anal.
Appl. \vol 88 \yr 1982 \pages 319--332 \endref
\ref \key \bf {\Bb} \by \bysame \paper 
Asymptotic behavior of nonexpansive mappings
\jour Contemp. Math. \vol 18 \yr 1983 
\pages 1--47 \endref
\ref \key  \bf{\DKR}  \by J. Dye, M. A. Khamsi and S.
Reich \paper  Random products of contractions in Banach spaces
\jour Trans.
Amer. Math. Soc. \vol 325\yr 1991 \pages 87--99\endref
%\ref \key \bf {\K} \by U.  Krengel \book Ergodic Theorems
%\publ De Gruyter \publaddr Berlin \yr 1985 \endref
\ref \key \bf {\DRa} \by J. Dye and S. Reich \paper Unrestricted
iterations of projections in Hilbert space \jour J. Math. Anal.
Appl. \vol 156 \yr 1991 \pages 101--119 \endref
\ref \key  \bf{\DRb}  \by \bysame
\paper Unrestricted iterations of nonexpansive mappings
in Hilbert space \jour  Nonlinear Analysis \vol 18\yr 1992 
\pages 199--207\endref
\ref \key \bf{\DRc}  \bysame \paper Unrestricted
iterations of nonexpansive mappings in Banach spaces
\jour Nonlinear Analysis \vol 18 \yr 1992 \pages 983--992
\endref

\ref \key  \bf{\vN}  \by J. von Neumann  \paper On
rings of operators.  Reduction theory \jour Ann. of Math. \vol
50\yr 1949 \pages 401--485 \endref 

\endRefs
\enddocument